\documentclass{article}
\usepackage{amssymb}
\newtheorem{thm}{Theorem}
\newtheorem{cor}[thm]{Corollary}
\newtheorem{lem}[thm]{Lemma}
\newtheorem{prp}[thm]{Proposition}
\newtheorem{???}[thm]{Conjecture}
\def\rem{\refstepcounter{thm}\noindent{\bf Remark \thethm\hskip1.7ex}}
\def\bull{\leavevmode\kern.1ex\vrule height1exwidth.9exdepth-.1ex\kern.8ex} 

\def\askip{\vskip\belowdisplayskip}
\def\block#1#2{{\setbox0=\hbox{#1\kern1ex}
           \leftskip=\wd0\parindent=-\wd0\par\leavevmode\box0 #2
           \par}}

\title{The size of bipartite graphs with girth eight}
\author{Stefan Neuwirth\footnotemark}

\begin{document}
\maketitle
\footnotetext[1]{
\noindent Centro Vito Volterra, Universit\`a Roma 2, Via di Tor
Vergata s.n.c., I-00133 Roma}

\begin{abstract}
  \noindent Reiman's inequality for the size of bipartite graphs of girth six is
  generalized to girth eight. It is optimal in as far as
  it admits the algebraic structure of generalized quadrangles as case
  of equality. This enables us to obtain the optimal estimate $e\sim
  v^{4/3}$ for balanced bipartite graphs. We also get an optimal
  estimate for very unbalanced graphs.
\end{abstract}

\section{Introduction}

De Caen and Sz\'ekely recently proposed a new bound for the size of a 
bipartite graph of girth eight, that is a bipartite graph without cycle of
length four and six. We adapt their method to obtain the following cubic
inequality.
\begin{thm}\label{TTHHMM}
  Let $G$ be a bipartite graph on $v+w$ vertices. 

\noindent $(i)$ If $G$ contains no
  cycle of length $4$ and $6$, then its size $e$ satisfies
$$    e^3-(v+w)e^2+2vwe-v^2w^2\le 0.$$
$(ii)$ If $v\ge\lfloor w^2/4\rfloor$, then furthermore $e\le v+\lfloor w^2/4\rfloor$.
\end{thm}
Part $(i)$ is the right generalization of Reiman's inequality for
bipartite graphs of girth $6$ (see Prop.~\ref{prp6}) to girth $8$. It
is optimal in the sense that it is an equality for all known extremal
graphs constructed {\it via\/} finite fields. Part $(ii)$ describes the case
of very unbalanced bipartite graphs and is optimal: there is a graph,
constructed by hand, for which it is an equality.

Let us give a brief description of this article. 
Section 2 describes a way to translate uncoloured graphs into
bipartite graphs and its converse. This permits to get two
propositions on very unbalanced graphs.

Section 3 summarizes facts about bipartite graphs of girth six that
should be folklore and well known although I did not see them printed.

Section 4 is the core of the paper. We adapt an inequality of
Atkinson {\it et al.\/} to get an optimal lower bound on the number of
paths of length $3$ in a bipartite graph
(Cor.~\ref{path3}). This enables us to bypass the final step in the proof
of \cite[Th.~1]{ds97} and to get our theorem.

Section 5 exploits the obtained results to get some handy information.

\section{Uncoloured graphs and bipartite graphs}
\label{sec:c}

\subsection{Expanding a graph to a bipartite graph}
\label{sec:c1}

We propose the following construction of a bipartite graph out of an
uncoloured graph. Let $G'$ be an uncoloured graph with set of vertices
$V$. Then the bipartite graph $G$ is defined as follows: \askip

\block\bull{the first class of vertices of $G$ is $V$;}

\block\bull{the second class $W$ of vertices of $G$ is the set of
  edges of $G'$;} 

\block\bull{the set of edges of $G$ is $\bigl\{\{x,y\}:y\hbox{ is an
    edge of $G'$ with endpoint }x\bigr\}.$}\askip

Thus every vertex of $W$ has degree $2$ and the size
of $G$ is twice the size of $G'$.

\subsection{Contracting a bipartite graph to an uncoloured graph}
\label{sec:c2}

Let us describe an inverse construction. 
Let $G$ be a bipartite graph with colour classes $V$ and $W$. 
Let $G'$ be the following graph: \askip

\block\bull{its set of vertices is $V$;}
\block\bull{its set of edges is 
$\bigl\{\{x,z\}\subseteq V:\exists
y~\{x,y\}{\rm~and~}\{z,y\}\hbox{ are edges of }G\bigr\}.$} 
\askip

The size of $G'$ is at most half the size of $G$. 
If $G$ contains no cycle of length $4$, then, given $\{x,z\}$, there
is at most one $y$ such that $\{x,y\}$ and $\{z,y\}$ are edges of $G$,
so that the size of $G'$ is exactly
\begin{equation}\label{sd}
\sum_{y\in W}{d(y)\choose 2}\le {\#V\choose2}.
\end{equation}
(We
recognize here \cite[Inequality $(2)$, p.~310]{bo78} for $s=t=2$.) 
Thus each vertex $y\in W$ of degree at least $2$ contributes at
least $1$ to sum $(\ref{sd})$. This yields
\begin{prp}\label{d2-4}
  Let $G$ be a bipartite graph on $v+w$ vertices that contains no
  cycle of length $4$. \vskip3pt

  \block{$\hphantom{i}(i)$}{If $w>{v\choose2}$, then there are at least
  $w-{v\choose2}$ vertices in $W$ of degree $0$ or $1$. }
  \block{$(ii)$}{If its minimal degree is at least $2$, then
  $w\le{v\choose2}$ and $v\le{w\choose2}$.}

\end{prp}
If $G$ contains no cycle of length $4$ nor $6$, then $G'$ contains
no triangle and its size is at most $\lfloor v^2/4\rfloor$.
This argument proves
\begin{prp}\label{d2-6}
  Let $G$ be a bipartite graph on $v+w$ vertices that contains no
  cycle of length $4$ or $6$. \vskip3pt

  \block{$\hphantom{i}(i)$}{If $w>\lfloor v^2/4\rfloor $, then there are at least
  $\lceil w-v^2/4\rceil$ vertices in $W$ with degree $0$ or $1$.}
  \block{$(ii)$}{If its minimal degree is at least $2$, then
  $w\le\lfloor v^2/4\rfloor$ and $v\le\lfloor w^2/4\rfloor$.}

\end{prp}

\section{Bipartite graphs of girth six}
\label{sec:4}

The following estimate is well known as Reiman's inequality, but its cases of
equality were not written down explicitly. Reading the proof of
\cite[Th.~VI.2.6]{bo78}, one gets with \cite[Def.~I.3.1]{bjl99}
\begin{prp}\label{prp6}
Let $v\le w$. A graph of girth at least $6$ on $v+w$ vertices with $e$
edges satisfies
$$O(v,w,e)=e^2-we-vw(v-1)\le0$$
$$e\le\sqrt{vw(v-1)+w^2/4}+w/2.$$
We have equality if and only if it is the incidence graph of a
Steiner system $S(2,k;v)$ on $v$ points with block degree $k$ given by $wk(k-1)=v(v-1)$.
\end{prp}
Note that by symmetry, we also get $O(w,v,e)\le0$, but this is
superfluous by 
\begin{lem}
  Let $v\le w$. Let $e$ be the positive root of $X^2-vX-vw(w-1)$. Then
  $O(v,w,e)\ge0$.
\end{lem}
{\it Proof. } As $(vw)^2-vvw-vw(w-1)=vw(vw-v-w+1)\ge0$,
we have $e\le vw$. Therefore
\begin{eqnarray*}
e^2-we-vw(v-1)&=&e^2-ve-vw(w-1)+(v-w)e+vw(w-v)\\&=&(vw-e)(w-v)\ge0
\end{eqnarray*}
\rem The case of equality in Prop.~\ref{prp6} implies the following:
By \cite[Cor.~I.2.11]{bjl99}, every vertex in $V$ has same degree $r$
and every vertex in $W$ has same degree $k$ with
\begin{equation}
  \label{eq:bl}
  k-1\mid v-1\quad{\rm and}\quad k(k-1)\mid v(v-1), 
\end{equation}
so that $v=1+r(k-1)$ and $k\mid r(r-1)$. For given $k$, this set of conditions is in
fact sufficient for the existence of an extremal graph for large
$r$: this is Wilson's Theorem \cite[Th.~XI.3.8]{bjl99}. For example,
we have the following complete sets of parameters
$(v,w,r,k)$: 
$$
(1+r(k-1),r(1+r(k-1))/k,r,k)\quad{\rm for}~1\le k\le 5~{\rm
  and}~k\mid r(r-1).$$
The first set of parameters satisfying $(\ref{eq:bl})$ for which an
extremal graph does not exist is $(36,42,7,6)$. Consult \cite[Table A1.1]{bjl99} for all known
block designs with $r\le 17$. \cite[Table A5.1]{bjl99} provides the
following sets of parameters $(v,w,r,k)$ for block designs: given any
prime power $q$ and natural number $n$, given $t\le s$,
$$\Bigl(q^n,q^{n-1}{q^n-1\over q-1},{q^n-1\over q-1},q\Bigr)\ ,\ 
\Bigl({q^{n+1}-1\over q-1},{q^{n+1}-1\over q^2-1}{q^n-1\over
  q-1},{q^n-1\over q-1},q+1\Bigr)\ ,$$
$$\bigl(q^3+1,q^2(q^2-q+1),q^2,q+1\bigr)\ ,\ 
\bigl(2^{t+s}-2^s+2^t,(2^s+1)(2^s-2^{s-t}+1),2^s+1,2^t\bigr)\ .$$ 

The following proposition provides a simpler but coarser bound.
\begin{prp}
  Let $G$ be a bipartite graph on vertex classes $V$ and $W$ with
  $\#V=v$ and $\#W=w$ without cycles of length $4$. Its size satisfies
$$e\le\cases{
\sqrt{2vw(v-1)}&if $w\le v(v-1)/2$\cr
v(v-1)/2+w&otherwise.\cr}$$
We have optimality in the second alternative for the
bipartite expansion of a complete graph on $V$ as described in Section
\ref{sec:c1}, on which we add
$w-v(v-1)/2$ new edges by connecting any vertex of $V$ to $w-v(v-1)/2$ new vertices in colour class $W$.
\end{prp}
{\it Proof. } 
By Proposition \ref{d2-4}, if $w>v(v-1)/2$, then $w-v(v-1)/2$ vertices in
$W$ have degree $0$ or $1$. If we remove them, we remove at most
$w-v(v-1)/2$ edges and the remaining graph has at most $v(v-1)$ edges because
$O(v,v(v-1)/2,v(v-1))=0$. The first alternative follows from
$$O(v,w,\sqrt{2vw(v-1)})=w\sqrt{v(v-1)}(\sqrt{v(v-1)}-\sqrt{2w}).\eqno{\bull}$$

\section{Bipartite graphs of girth eight}
\label{sec:6}

\subsection{Statement of the theorem}
\label{sec:ss1}

Consult \cite[Def.~1.3,1]{va98} for the definition of weak generalized
polygons.
\begin{thm}\label{thm:1}
  Let $G$ be a bipartite graph on vertex classes $V$ and $W$ with
  $\#V=v$ and $\#W=w$. If $G$ contains no cycle of length $4$ or $6$,
  then its size $e$ satisfies
  \begin{equation}
    \label{pi}
    P(v,w,e)=e^3-(v+w)e^2+2vwe-v^2w^2\le 0.
  \end{equation}
  We have equality exactly in two cases:

  \block{$\hphantom{i}(i)$}{if $G$ is the complete bipartite graph and $v=1$ or $w=1$.}
  
  \block{$(ii)$}{if $G$ is the incidence graph of a weak generalized quadrangle.}
\end{thm}

\rem Let us first note that this polynomial has exactly one positive
root in $e$ for positive $v,w$. It suffices to this purpose to show
that its discriminant is negative. This is $-v^2w^2D$ with
$$D=27p^2+4s^3-36sp-4s^2+32p\quad,\quad s=v+w,\,p=vw.$$
Let us study this quantity for $s\ge2$, $p\ge s-1$. We have
$$
{dD\over dp}=54p-36s+32\ge54p-36(p+1)+32=18p-4>0,
$$
so that its minimum satisfies $p=s-1$, which implies
$D=(4s-5)(s-1)^2\ge3$. Therefore Inequality $(\ref{pi})$ is equivalent
to an inequality of form $e\le e(v,w)$.\askip

\rem The case of equality in Th.~\ref{thm:1} implies the following:
every vertex in $V$ has same degree $s+1$ and every vertex in $W$ has
same degree $t+1$. By \cite[Cor.~1.5.5, Th.~1.7.1]{va98}, $s+t\mid
st(1+st)$ and
$$v=(t+1)(1+st)~,~w=(s+1)(1+st)~,~e=(s+1)(t+1)(1+st).$$
Let us suppose, by symmetry, that $s\le t$. If $s=0$, we get case $(i)$. If $s=1$, we obtain exactly the
examples of extremal graphs produced by de Caen and Sz{\'e}kely:
$W$ consists of $t+1$ horizontal lines and as much vertical
lines and $V$ is the set of $(t+1)^2$ intersection points and
$G$ is the point-line incidence graph of this grid (this is also
the bipartite expansion of a complete bipartite graph on
$(t+1)+(t+1)$ vertices.)  Otherwise $s,t\ge2$ and $G$ is
in fact the incidence graph of a generalized quadrangle, so that
by \cite[Th.~1.7.2]{va98}, $t\le s^2$. Let
$q$ be a prime power. Then there are generalized quadrangles with
set of parameters $(s,t)$ any of $(q,q)$, $(q,q^2)$, $(q^2,q^3)$,
$(q-1,q+1)$; all known ones fit in this list. In
particular, by
\cite[Th.~1.7.9]{va98}, if $t\ge s=2$, then $t=2$ or
$t=4$ and in each case there is exactly one extremal graph. By
\cite[Sec.~1.7.11]{va98}, if $t\ge s=3$, then there is a (unique)
extremal graph exactly if $t=3,5,9$. There is a unique extremal
graph with $s=t=4$. It is open whether there exists a
generalized quadrangle with $s=4$ and $t\in\{11,12\}$.  \askip

\subsection{A generalization of an inequality of Atkinson {\it et al.}}
\label{sec:ss2}

We first need an optimal lower bound on the number of paths of length
3.  Let us prove the following inequality.

\begin{thm}\label{AWM}
  Let $(a_{ij})_{1\le i\le v,1\le j\le w}$ be a matrix of nonnegative
  coefficients and $\rho,\gamma\ge0$. Let 
\begin{equation}
  \label{eq:part}
  a_{i\star }=\sum_{j=1}^wa_{ij}\quad,\quad 
  a_{\star j}=\sum_{i=1}^va_{ij}\quad, \quad e=\sum_{i=1}^n\sum_{j=1}^va_{ij}.
\end{equation}
  If $a_{i\star }\ge2\rho$ and
  $a_{\star j}\ge2\gamma$, then
  \begin{equation}
    \label{eq:awmbis}
    \phi=\sum_{i=1}^v\sum_{j=1}^wa_{ij}
  (a_{i\star}-\rho)(a_{\star j}-\gamma)\ge 
e(e/v-\rho)(e/w-\gamma),
  \end{equation}
  equality holding exactly if $a_{i\star }$ and $a_{\star j}$ are constant.
\end{thm}
This refines the inequality in \cite{awm60}, which states
\begin{equation}
  \label{eq:awm}
  \psi=\sum_{i=1}^v\sum_{j=1}^wa_{ij}a_{i\star}a_{\star j}\ge {e^3/ vw}
\end{equation}
as, by the Arithmetic-Quadratic Mean Inequality,
\begin{equation}
  \label{eq:refin}
  \phi-\psi=-\gamma\sum_{i=1}^va_{i\star}^2-\rho\sum_{j=1}^wa_{\star
  j}^2+\rho\gamma e\le e(-\gamma e/v-\rho e/w+\rho\gamma)
\end{equation}

\rem 
If $v=w$ and $a$ is diagonal, Inequality $(\ref{eq:awm})$ is the
Arithmetic-Cubic Mean Inequality and Inequality $(\ref{eq:awmbis})$
becomes
$$
{1\over v}\sum_{i=1}^va_{ii}(a_{ii}-\rho)(a_{ii}-\gamma)\ge{e\over
  v}{e-v\rho\over v}{e-v\gamma\over v},$$
which is true by Tchebychef's Inequality \cite[Th.~43]{hlp52} if $a_{ii}\ge\rho$ and
$a_{ii}\ge\gamma$. For our ``non commutative Tchebychef Inequality'', the conditions $a_{i\star }\ge2\rho$ and $a_{\star j}\ge2\gamma$
cannot be weakened to $a_{i\star }\ge\rho$ and $a_{\star j}\ge\gamma$,
as we have the following counterexamples:  
$$\pmatrix{2&5\cr4&0\cr},\pmatrix{0&1&1\cr1&0&0\cr1&0&0\cr}.$$
\askip

{\it Proof. } If $(\ref{eq:awmbis})$ is an equality, then so are
$(\ref{eq:refin})$ and $(\ref{eq:awm})$ and our case of equality
follows from the identical case of equality in \cite{awm60}, whose proof we now imitate. We shall suppose
that $a_{i\star }>2\rho$ or $a_{\star j}>2\gamma$, so that the
whole inequality follows by continuity. Fix $e$ and suppose that under
this condition the $a_{ij}$ are chosen so to minimize $\phi$. We may
suppose that the rows and the columns have been permuted such that the
sequences $(a_{i\star})$ and $(a_{\star j})$ are nondecreasing:
\begin{equation}\label{dec}
a_{1\star}\le\dots\le a_{v\star}\quad,\quad 
a_{\star 1}\le\dots\le a_{\star w}.
\end{equation}
If one of these sequences is constant, the inequality follows by the
Arithmetic-Quadratic inequality (and the case of equality is easy).
Let us suppose that this is not so.

One can suppose that $a_{1w}$ and $a_{v1}$ are positive.  Let us show
the argument for $a_{1w}$. If $a_{1w}=0$, there are $k,l$ such that
$a_{1k},a_{lw}>0$. Make a perturbation by adding $\alpha$ to $a_{1w}$
and to $a_{lk}$ and subtracting $\alpha$ to $a_{1k}$ and to $a_{lw}$.
The row and column sums $a_{i\star }$ and $a_{\star j}$ are unaltered
and $\phi$ increases of
\begin{eqnarray*}
\Delta\phi&=&\alpha\bigl(
   (a_{1\star}-\rho)(a_{\star w}-\gamma)+(a_{l\star}-\rho)(a_{\star k}-\gamma)\\
&&\mathop{-}(a_{1\star}-\rho)(a_{\star k}-\gamma)-(a_{l\star}-\rho)(a_{\star m}-\gamma)\bigr)\\
&=&\alpha
(a_{1\star}a_{\star w}+a_{l\star}a_{\star k}-a_{1\star }a_{\star k}-a_{l\star }a_{\star w})\\
&=&\alpha(a_{1\star}-a_{l\star})(a_{\star w}-a_{\star k}),
\end{eqnarray*}
so that $\phi$ does not increase.

Now make the following perturbation: add $2\alpha$ to $a_{11}$ and subtract
$\alpha$ to $a_{1w}$ and to $a_{v1}$. Let us compute the differential
of $\phi$: as 
$${d\phi\over da_{rc}} = (a_{r\star }-\rho)(a_{\star c}-\gamma) +
\sum_{i=1}^va_{ic}(a_{i\star }-\rho) + \sum_{j=1}^wa_{rj}(a_{\star j}-\gamma),$$
\begin{eqnarray*}
d\phi&=&d\alpha\Bigl(2{d\phi\over da_{11}}-{d\phi\over da_{1w}}-{d\phi\over da_{v1}}\Bigr)\\
&=&d\alpha\Bigl((a_{1\star }-\rho)(a_{\star 1}-a_{\star w}) +
(a_{1\star }-a_{v\star })(a_{\star 1}-\gamma) +
\sum_{i=1}^va_{i1}(a_{i\star }-\rho) 
\\&&\qquad\qquad\mathop{+}\sum_{j=1}^wa_{1j}(a_{\star j}-\gamma) - 
\sum_{i=1}^va_{iw}(a_{i\star }-\rho) - 
\sum_{j=1}^wa_{vj}(a_{\star j}-\gamma)\Bigr)\\
\end{eqnarray*}
For positive $d\alpha$, we have by $(\ref{dec})$
\begin{eqnarray*}
d\phi&\le&d\alpha\bigl(
(a_{1\star }-\rho)(a_{\star 1}-a_{\star w}) +
(a_{1\star }-a_{v\star })(a_{\star 1}-\gamma) +
a_{\star 1}(a_{v\star }-\rho)\\
&&\qquad\qquad\qquad\quad
\mathop{+}a_{1\star }(a_{\star w}-\gamma)-a_{\star w}(a_{1\star }-\rho)
-a_{v\star }(a_{\star 1}-\gamma)\bigr)\\
&=&d\alpha\bigl(
(a_{1\star }-2\rho)(a_{\star 1}-a_{\star w})+(a_{\star
  1}-2\gamma)(a_{1\star }-a_{v\star })\bigr)\\
&<&0\,,
\end{eqnarray*}
which contradicts the minimum hypothesis.\hfill$\bull$\askip

\begin{cor}\label{path3}
  Let $G$ be a bipartite graph on $v+w$ vertices and of minimal degree
  $2$. Then the number of paths of length $3$ in $G$ is at least
  $e(e/v-1)(e/w-1)$. This bound is achieved exactly if the graph is
  regular for each of its two colours.
\end{cor}
{\it Proof. } 
A path of length 3 is a sequence of 4 vertices
$(x,y,z,t)$ with no repetition such that 
$$\{x,y\},\{y,z\},\{z,t\}\in G.$$ 
Given two adjacent vertices $y$ and $z$, the number of
paths $(x,y,z,t)$ makes $(d(y)-1)(d(z)-1)$, where $d$ denotes the
degree of a vertex. Therefore the number of all paths of length
$3$ is $$\sum_{\{y,z\}\in G}(d(y)-1)(d(z)-1).$$ 
Let $(a_{ij})_{1\le i\le v,1\le j\le w}$ be the reduced incidence
matrix of $G$: $a_{ij}=1$ if the $i$th vertex of the first class is
adjacent to the $j$th vertex of the second class; otherwise $a_{ij}=0$. Then this sum is 
\begin{equation}
  \label{eq:sum}
  \sum_{i=1}^v\sum_{j=1}^wa_{ij}(a_{i\star }-1)(a_{\star j}-1),
\end{equation}
so that it suffices to take $\rho=\gamma=1$ in Th.~\ref{AWM}.

\subsection{Proof of Theorem \ref{thm:1}}
\label{sec:proof}

The case of equality follows from \cite[Lemma~1.4.1]{va98}
because its axiom $(i)$ is exactly what makes Bound
$(\ref{eq:bound})$ an equality.

I now follow the proof of \cite[Th.~1]{ds97}. If $G$ contains no cycle of
length $4$, there is no path of length $3$ between two adjacent
vertices; if $G$ contain no cycle of
length $6$, there is at most one path of length $3$ between non-adjacent
vertices of different colour. Therefore the sum $(\ref{eq:sum})$ is bounded by 
\begin{equation}
  \label{eq:bound}
  vw-e\quad{\rm with}\quad e=\sum_{i=1}^n\sum_{j=1}^va_{ij}.
\end{equation}
By Corollary~\ref{path3}, if all the
vertices of $G$ have degree at least two, one has 
$$vw-e\ge e^3/vw-(1/v+1/w)e^2+e$$
and therefore $(\ref{pi})$.  In
order to get rid of this degree condition, we have to do an induction
on the sum $s=v+w$ of the number of columns and the number of rows of
the incidence matrix. If $v=1$, then $P(v,w,e)=(e-w)(e^2-e+w)$, so
that the inequality states $e\le w$, which is trivial; symmetrically
for $w=1$. Suppose the result is true for all $v\times w$ incidence
matrices with $v+w=s$.  Consider now a $v\times w$ incidence matrix
with $v+w=s+1$ and $v,w\ge2$. If each vertex has degree at least two,
the result is true; otherwise there is a column or a row containing
only zeroes or exactly one ``1''. Apply the induction hypothesis on the matrix
without this row or column: we get $P(v-1,w,e-1)\le0$ or
$P(v,w-1,e-1)\le0$ and we may apply the following growth lemma to
conclude.\hfill\bull

\begin{lem}
Let $v,w\ge1$. If $P(v,w,e)\le0$, then $P(v+1,w,e+1)\le0$.
\end{lem}
In fact, one has
$$P(v+1,w,e+1)-P(v,w,e)=2{e}^{2}+(1-2v)e+(w-{w}^{2})(2v+1)-v,$$ which is negative as
long as 
$$0\le e\le e_0=\bigr(2v-1+\sqrt{ \left (2v+1\right )\left
    (2v+8{w}^{2}-8w+1\right )}\bigl)/4=(2v-1+\Delta)/4.$$
Let us use that $P(v,w,e)$ has a unique root in $e$ and compute
$P(v,w,e_0)$. This makes
$$(4vw^2+2w^2+1)\Delta/16+(-16vw^3-8v^2w^2-8w^3+8vw^2+2w^2-2v+4w-1)/16.$$
Then either the second term in this sum is positive and $P(v,w,e_0)$
is positive, or the conjugate expression of this sum is positive, and
the product of the sum with this conjugate expression is
$$(w-1)^2w^2(8v^3w^2+4v^2w^2-2vw^2+2v^2-w^2-4vw+2v-2w)/8,$$ which is
positive if $v,w\ge1$.\hfill\bull

\subsection{Further remarks}
\label{sec:ss3}

\rem 
Theorem \ref{thm:1} does not always give the right order of
magnitude for the maximal size of a graph of girth $8$: as 
$$P(v,w,(vw)^{2/3})=2(vw)^{5/3}-(vw)^{4/3}(v+w))\le2(vw)^{5/3}-2(vw)^{4/3+1/2}\le0,$$
we expect to find maximal graphs of size $(vw)^{2/3}$: De~Caen and
Sz\'ekely \cite[Th.~4]{ds97} find a counterexample to this expectation if $v$ ``lies in
an interval just slightly below'' $w^2$. They conjecture \cite{ds92}
that this is the case as soon as
$v\gg w^{5/4}$ and $v\ll w^2$.\askip

In the case of $v=w$, let us give the following approximation for the
real root of the polynomial. For
$$e=v^{4/3}+{2\over3}v-{2\over9}v^{2/3}-{20\over81}v^{1/3},$$
$$P(v,v,e)={\frac {40}{243}}{v}^{7/3}+{\frac
  {376}{2187}}{v}^{2}-{\frac {80}{2187}}{v}^{5/3}-{\frac
  {800}{19683}}{v}^{4/3}-{\frac {8000}{531441}}v\ge{\frac
  {129808}{531441}},$$
\begin{eqnarray*}
P(v,v,e-16/81)&=&-{8\over531441}(v^{1/3}-1)\bigl(
39366{v}^{7/3}+28431{v}^{2}+8262{v}^{5/3}\\
&&\mathop{-}8748{v}^{
4/3}-11880{v}-6560{v}^{2/3}-2432v^{1/3}-512\bigr)\\&\le&0.
\end{eqnarray*}
In particular,
\begin{cor}
  Let $G$ be a bipartite graph of size $e$ with $v$ vertices in each vertex
  class. If the girth of $G$ is at least $8$, then
$$e<v^{4/3}+{2\over3}v-{2\over9}v^{2/3}-{20\over81}v^{1/3}.$$
\end{cor}

Let us now show that we generalize the following estimations for the
size of bipartite graphs of girth $8$ in \cite[Th.~1]{ds97}:\askip

\block{$\hphantom{i}(i)$}{\it if the minimal degree of $G$ is at least $2$, then
  $e\le2^{1/3}(vw)^{2/3}$;} 

\block{$(ii)$}{\it if $v\preccurlyeq w^2$ or $w\preccurlyeq v^2$, then
$e\preccurlyeq(vw)^{2/3}$.}\askip

\noindent In fact,
$$P(v,w,2^{1/3}(vw)^{2/3})=(vw)^{4/3}(w^{2/3}-2^{2/3}v^{1/3})(v^{2/3}-2^{2/3}w^{1/3}),$$
which is nonnegative exactly if $v\le w^2/4$ and $w\le v^2/4$ or if
$(v,w)$ is among $\bigl\{(1,1),(1,2),(2,1),(2,2),(3,3)\bigr\}$, and
this is the case by Prop.~\ref{d2-6} if the minimal degree is at least
$2$.

Furthermore, by Prop.~\ref{d2-6}, if $w>\lfloor v^2/4\rfloor$, then
$\lceil w-v^2/4\rceil$ vertices in
$W$ have degree $0$ or $1$. If we remove them, we remove at most
$\lceil w-v^2/4\rceil$ edges and the remaining graph has at most
$\lfloor v^2/2\rfloor$ edges because $P(v,\lfloor v^2/4\rfloor,\lfloor v^2/2\rfloor+1)>0$. 
This yields
\begin{prp}
  Let $G$ be a bipartite graph on vertex classes $V$ and $W$ with
  $\#V=v$ and $\#W=w$ without cycles of length $4$ and $6$. Its size satisfies
$$e\le\cases{ 2^{1/3}(vw)^{2/3}&if $\max(v,w)\le \lfloor
  \min(v,w)^2/4\rfloor$\cr \lfloor
  \min(v,w)^2/4\rfloor+\max(v,w)&otherwise.\cr}$$
We have optimality in the second alternative: make
a bipartition $V=V_1\cup V_2$ with $V_1=\lceil v/2\rceil$ and
$V_2=\lfloor v/2\rfloor$, let $G'$ be the complete bipartite graph
on the colour classes $V_1$ and $V_2$, which has $\lfloor
v^2/4\rfloor$ edges. Now consider the bipartite expansion of $G'$,
add $\lceil w-v^2/4\rceil$ new vertices to colour class $W$,
and connect each of them to some vertex of $V$.
\end{prp}
Note that this estimate yields another proof of \cite[Th.~1]{gy97} by
means of \cite[Th.~3]{gy97}.  \askip


\rem Our inequality condenses the following facts about the behaviour of $e$ for
fixed $w$ and large $v$. If $w\le 3$, then extremal graphs of girth
$8$ do not contain any cycle at all, so that their size is $e=v+w-1$;
if $v\ge w=4$ and if $v=w=5$, then extremal graphs of girth
$8$ contain exactly one cycle, so that their size is $e=v+w$; if
$v>w=5$, then extremal graphs of girth $8$ contain exactly one ``$\theta$-graph'',
so that their size is $e=v+w+1$.


\end{document}